\documentclass[11pt]{article}
\usepackage{mathrsfs}
\usepackage{amscd}
\usepackage{amsmath,amsfonts,amssymb,amscd}
\usepackage{indentfirst,graphics,epsfig,psfrag}
\input{epsf}
\usepackage{ifpdf}
\usepackage{enumerate}
\usepackage{appendix}
\usepackage{enumerate}
\usepackage{color}
\usepackage{lineno}

\usepackage{lineno}
\makeatletter \@addtoreset{figure}{section} \makeatother
\makeatletter
\long\def\@makecaption#1#2{%
   \vskip 10\p@
   \setbox\@tempboxa\hbox{{#1}\ \ #2}%
   \ifdim \wd\@tempboxa >\hsize

       {#1}\ \ #2\par
   \else
       \hbox to\hsize{\hfil\box\@tempboxa\hfil}%
   \fi}
\makeatother

\newtheorem{thm}{Theorem}
\newtheorem{cor}{Corollary}
\newtheorem{lem}{Lemma}

\newtheorem{pro}{Proposition}

\newcommand{\qed}{{\hfill\rule{3pt}{7pt}}}

\def\qed{\hfill \rule{4pt}{7pt}}

\begin{document}
\title{\textbf{Proper connection number and graph products}
\footnote{Supported by the National Science Foundation of China (No.
11161037) and the Science Found of Qinghai Province (No.
2014-ZJ-907).}}
\author{
\small Yaping Mao, \  Fengnan Yanling \footnote{Corresponding author},  \  Zhao Wang,  \  Chengfu Ye\\[0.3cm]
\small Department of Mathematics, Qinghai Normal\\
\small University, Xining, Qinghai 810008, China}
\date{}
\maketitle
\begin{abstract}
A path $P$ in an edge-colored graph $G$ is called \emph{a proper path} if no two adjacent
edges of $P$ are colored the same, and $G$ is \emph{proper connected} if every two vertices
of $G$ are connected by a proper path in $G$. The \emph{proper connection number} of a
connected graph $G$, denoted by $pc(G)$, is the minimum number of colors that are
needed to make $G$ proper connected. In this paper, we study the proper connection number on the lexicographical, strong, Cartesian,
and direct product and present several upper bounds for these
products of graphs.
\\[2mm]
\end{abstract}

\section{Introduction}

All graphs considered in this paper are simple, finite and undirected. We follow the
terminology and notation of Bondy and Murty \cite{Bondy}. For a graph $G$, we use $V(G)$, $E(G)$,
$n(G)$, $m(G)$, $\delta(G)$, $\kappa(G)$, $\kappa'(G)$, $\delta(G)$ and $diam(G)$ to denote the vertex set, edge set, number of vertices, number of edges, connectivity, edge-connectivity, minimum degree and diameter of $G$, respectively. The rainbow connections of a graph which are applied to measure the safety of a network are introduced by Chartrand, Johns, McKeon and
Zhang \cite{Chartrand}. Readers can see \cite{Chartrand, ChartrandJMZ,
ChartrandOZ} for details. An edge-coloring of a graph $G$ is an assignment $c$ of colors to the edges of $G$, one color
to each edge of $G$. Consider an edge-coloring (not necessarily
proper) of a graph $G=(V,E)$. We say that a path of $G$ is
\emph{rainbow}, if no two edges on the path have the same color. An
edge-colored graph $G$ is \emph{rainbow connected} if every two
vertices are connected by a rainbow path. The minimum number of
colors required to rainbow color a graph $G$ is called \emph{the
rainbow connection number}, denoted by $rc(G)$. For more
results on the rainbow connection, we
refer to the survey paper \cite{LiSS} of Li, Shi and Sun and a new
book \cite{LiS} of Li and Sun.

If adjacent edges of $G$ are assigned different colors by $c$, then $c$ is a
proper (edge-)coloring. The minimum number of colors needed in a proper coloring of G is
referred to as the chromatic index of G and denoted by $\chi'(G)$. Recently, Andrews, Laforge, Lumduanhom and Zhang \cite{AndrewsLLZ}
introduce the concept of proper-path colorings. Let $G$ be an edge-colored graph, where
adjacent edges may be colored the same. A path $P$ in $G$ is called a \emph{proper path} if no two
adjacent edges of $P$ are colored the same. An edge-coloring $c$ is a \emph{proper-path coloring} of
a connected graph $G$ if every pair of distinct vertices $u,v$ of $G$ is connected by a proper
$u$-$v$ path in $G$. A graph with a proper-path coloring is said to be \emph{proper connected}. If $k$
colors are used, then $c$ is referred to as a proper-path $k$-coloring. The minimum number
of colors needed to produce a proper-path coloring of $G$ is called the \emph{proper connection
number} of $G$, denoted by $pc(G)$.

Let $G$ be a nontrivial connected graph of order $n$ and size $m$. Then the proper
connection number of $G$ has the following bounds.
$$
1\leq pc(G)\leq \min\{\chi'(G),rc(G)\}\leq m.
$$
Furthermore, $pc(G)=1$ if and only if $G=K_n$ and $pc(G)=m$ if and only if $G = K_{1,m}$ is
a star of order $m+1$. For more details on the proper
connection number, we refer to \cite{AndrewsLLZ, LiWY, HuangLW}.

The standard products (Cartesian, direct, strong, and lexicographic)
draw a constant attention of graph research community, see some
recent papers \cite{AnandCKP, Klavzar, Nowakowski, Zhu}.

In this paper, we consider four standard products: the
lexicographic, the strong, the Cartesian and the direct with respect
to the proper connection number. Every of these
four products will be treated in one of the forthcoming sections.

\section{The Cartesian product}

The \emph{Cartesian product} of two graphs $G$ and $H$, written as
$G\Box H$, is the graph with vertex set $V(G)\times V(H)$, in which
two vertices $(g,h)$ and $(g',h')$ are adjacent if and only if
$g=g'$ and $(h,h')\in E(H)$, or $h=h'$ and $(g,g')\in E(G)$.
Clearly, the Cartesian product is commutative, that is, $G\Box H$ is
isomorphic to $H\Box G$.

\begin{lem}{\upshape \cite{Hammack}}\label{lem1}
Let $gh$ and $g'h'$ be two vertices of $G\Box H$. Then
$$
d_{G\Box H}(gh,g'h')=d_{G}(gg')+d_{H}(hh').
$$
\end{lem}

\begin{thm}\label{th1}
Let $G$ and $H$ be connected graphs with $|V(G)|\geq 2$ and $|V(H)|\geq 2$. Then
$$
pc(G\Box H)\leq \min\{pc(G), pc(H)\}+1.
$$
Moreover, the bound is sharp.
\end{thm}
\begin{pf}
Without loss of generality, we assume $pc(H)\leq pc(G)$.
Suppose $\{0,1,\cdots,pc(H)-1\}$ be a proper coloring
of $H$. Clearly, Since $G$ is connected, there is a path connecting $g$ and $g'$, say
$P=gg_{1},\ldots g_{\ell-1}g'$ where $g'=g_{\ell}$. By the same reason, there is a
path connecting $h$ and $h'$, say $Q=hh_{1},\ldots h_{k-1}h'$ where $h'=h_{k}$.
Now we give a coloring of $G\Box H$ using $pc(H)+1$ colors. To show that
$pc(G\Box H)\leq pc(H)+1$, we provide a proper-coloring $c$ of $G\Box H$ with $pc(H)+1$ colors as follows.
$$
\left\{
\begin{array}{ll}
c(gh_s,gh_t)=c(h_sh_t), &~if~s\neq t.\\[0.2cm]
c(g_ih,g_jh)=pc(H)+1, &~if~i\neq j;\\[0.2cm]
\end{array}
\right.
$$
It suffices to check that there is a proper-path between any two vertices
$(g,h),(g',h')$ in $G\Box H$. If $g=g'$ or $h=h'$,
then $P$ or $Q$, respectively, is a trivial one vertex path.
We distinguish the following two cases to prove this theorem.

{\bf Case $1$.} $h=h'$

If $\ell$ is even, then we let $h_1$ be an arbitrary neighbor of $h$.
The path induced by the edges in
$$
\{(gh, g_1h), (g_1h, g_1h_1), (g_1h_1, g_2h_1),\cdots, (g_{\ell-1}h_1, g'h_1), (g'h_1, g'h')\}
$$
is proper $(g, h), (g', h')$-path in $G\Box  H$.

If $\ell$ is odd, then we let $h_1$ be an arbitrary neighbor of $h$.
The path induced by the edges in
$$
\{(gh, g_1h), (g_1h, g_1h_1), (g_1h_1, g_2h_1),\cdots, (g_{\ell-1}h_1, g_{\ell-1}h), (g_{\ell-1}h, g'h')\}
$$
is proper $(g, h), (g', h')$-path in $G\Box  H$.

{\bf Case $2$.} $h\neq h'$

If $g=g'$, then $(g,h),(g',h')\in
H(g)$. Clearly, there is a proper-path connecting $(g,h)$ and $(g',h')$.
Now we consider $g\neq g'$.
If $\ell$ is even, then we let $h_1$ be an arbitrary neighbor of $h$.
The path induced by the edges in
$$
\{(gh, g_1h), (g_1h, g_1h_1), (g_1h_1, g_2h_1),\cdots, (g_{\ell-1}h_1, g_{\ell}h_1), (g'h_1, g'h_2)
$$,
$$
(g'h_2, g'h_3)\cdots (g'h_{k-1}, g'h')\}
$$

is proper $(g, h), (g', h')$-path in $G\Box  H$.

If $\ell$ is odd, then we let $h_1$ be an arbitrary neighbor of $h$.
The path induced by the edges in
$$
\{(gh, g_1h), (g_1h, g_1h_1), (g_1h_1, g_2h_1),\cdots, (g_{\ell-1}h_1, g_{\ell-1}h), (g_{\ell-1}h, g'h)
$$
$$
(g'h, g'h_1),\cdots, (g'h_{k-1}, g'h')\}
$$
is a proper $(g, h), (g', h')$-path in $G\Box  H$.\qed
\end{pf}\\

To show the sharpness of the above bound, we consider the following
example.\\

\noindent \textbf{Example 1:} Let $G=P_2$ and $H=K_n$. Then
$pc(G\Box H)\leq \min\{pc(G), pc(H)\}+1=2$ by Theorem \ref{th1}. From Lemma \ref{lem1}, we have $diam(G\Box H)=diam(G)+diam(H)=2$ and hence $pc(G\Box H)\geq 2$. Therefore,
$pc(G\Box H)=2=\min\{pc(G), pc(H)\}+1$.

\section{The strong product}

The \emph{strong product} $G\boxtimes H$ of graphs $G$ and $H$ has
the vertex set $V(G)\times V(H)$. Two vertices $(g,h)$ and $(g',h')$
are adjacent whenever $gg'\in E(G)$ and $h=h'$, or $g=g'$ and $hh'
\in E(H)$, or $gg'\in E(G)$ and $hh'\in E(H)$.

\begin{lem}{\upshape \cite{Hammack}}\label{lem2}
If $G$ is a nontrivial connected graph and $H$ is a connected spanning subgraph of $G$, then $pc(G)\leq pc(H)$.
\end{lem}

The strong product is connected whenever both factors are and the
vertex connectivity of the strong product was solved recently by
Spacapan in [23].

By Lemma \ref{lem2}, we have $pc(G\boxtimes H)\leq pc(G\Box H)$. By Theorem \ref{th1}, the following proposition is immediate.

\begin{pro}\label{pro1}
Let $G$ and $H$ be connected graphs. Then
$$
pc(G\boxtimes H)\leq \min\{pc(G), pc(H)\}+1.
$$
Moreover, the bound is sharp.
\end{pro}

\begin{lem}{\upshape \cite{Hammack}}\label{lem3}
Let $gh$ and $g'h'$ be two vertices of $G\Box H$. Then
$$
d_{G\boxtimes H}(gh,g'h')=\max\{d_{G}(gg'), d_{H}(hh')\}.
$$
\end{lem}

To show the sharpness of the upper bound in Proposition \ref{pro1}, we consider the following
example.\\

\noindent \textbf{Example 2:} Let $G=P_n$ be a complete graph and $H=P_2$.
From Proposition \ref{pro1}, we have $pc(G\boxtimes H)\leq \min\{pc(G), pc(H)\}+1=2$ .
By Lemma \ref{lem3}, $diam(G\boxtimes H)\geq 2$ and hence $pc(G\boxtimes H)\leq 2$.
Therefore, $pc(G\boxtimes H)=2= \min\{pc(G), pc(H)\}+1$.\\

\section{The lexicographical product}

The lexicographic product $G\circ H$ of graphs $G$ and $H$ has the
vertex set $V(G\circ H)=V(G)\times V(H)$. Two vertices
$(g,h),(g',h')$ are adjacent if $gg'\in E(G)$, or if $g=g'$ and
$hh'\in E(H)$. The lexicographic product is not commutative and is
connected whenever $G$ is connected.

In this section, let $G$ and $H$ be two connected graphs with
$V(G)=\{g_1,g_2,\ldots,g_{n}\}$ and $V(H)=\{h_1,h_2,\ldots,h_{m}\}$,
respectively. Then $V(G\circ H)=\{(g_i,h_j)\,|\,1\leq i\leq n, \
1\leq j\leq m\}$. For $h\in V(H)$, we use $G(h)$ to denote the
subgraph of $G\circ H$ induced by the vertex set
$\{(g_i,h)\,|\,1\leq i\leq n\}$. Similarly, for $g\in V(G)$, we use
$H(g)$ to denote the subgraph of $G\circ H$ induced by the vertex
set $\{(g,h_j)\,|\,1\leq j\leq m\}$.

\begin{thm}\label{th2}
Let $G$ and $H$ be connected graphs.

$(i)$ For $pc(G), pc(H)\geq 2$, we have
$$
\left\{
\begin{array}{ll}
pc(G\circ H)\leq pc(H), &~if~pc(G)>pc(H);\\[0.2cm]
pc(G\circ H)\leq pc(G)+1, &~if~pc(G)<pc(H);\\[0.2cm]
pc(G\circ H)\leq pc(G), &~if~pc(G)=pc(H).
\end{array}
\right.
$$

$(ii)$ If $pc(G)=1, pc(H)\geq 2$, then $pc(G\circ H)=2$;

$(iii)$ If $pc(H)=1, pc(G)\geq 2$, then $pc(G\circ H)=2$;

$(iv)$ If $pc(G)=1, pc(H)=1$, then $pc(G\circ H)=1$.

Moreover, the bound is sharp.
\end{thm}
\begin{pf}
$(i)$ If $pc(G)>pc(H)$, then we give a coloring of $G\circ H$ using $pc(H)$ colors.
Suppose $pc(H)=\{1,2,\cdots, pc(H)\}$ is a proper-coloring of $H$.
We color the edges $c(gh_i,gh_j)$ $(i\neq j)$ the same as $H$, and the edges $c(g_ih_s,g_jh_t)=1$ $(i\neq j)$.
It suffices to
check that there is a proper-path between any two vertices
$(g,h),(g',h')$ in $G\circ H$. If $g=g'$, then there is a proper path in $H(g)$ as desired. Now suppose $g\neq g'$. Since $pc(H)\geq 2$, there is an edge $h_ih_j\in E(H)$ such that $c(h_ih_j)\neq 1$. The path induced by the edges in
$$
\{(gh, gh_1), (gh_1, gh_2)\cdots (gh_{i-1}gh_{j-1}, gh_igh_j), (gh_igh_j, g_1h_j), (g_1h_j,g_1h_i),
$$
$$
(g_1h_i,g_2h_j), (g_2h_j,g_2h_i),\cdots (g_{\ell-1}h_j, g'h')\}
$$
is a proper-path connected $gh$ and $g'h'$.

If $pc(G)<pc(H)$, then $pc(G\circ H)\leq pc(G)+1$ by Lemma \ref{lem2} and Theorem \ref{th1}.

If $pc(G)=pc(H)$, then we color $G\circ H$ as follows.
$$
\left\{
\begin{array}{ll}
c(g_ih,g_jh)=c(g_ig_j), &~if~i\neq j;\\[0.2cm]
c(gh_s,gh_t)=c(h_sh_t), &~if~s\neq t;\\[0.2cm]
c(g_ih_s,g_jh_t)=c(g_ig_j), &~if~i\neq j\ and\ s\neq t.
\end{array}
\right.
$$
It suffices to
check that there is a proper-path between any two vertices
$(g,h),(g',h')$ in $G\circ H$.
If $h=h'$, then there is a
proper-path connecting $(g,h)$ and $(g',h')$ in $G(h)$, as
desired. Suppose $h\neq h'$. If $g=g'$, then $(g,h),(g',h')\in
H(g)$. There is a
proper-path connecting $(g,h)$ and $(g',h')$.
We now assume $g\neq g'$. Since $G$ is
connected, it follows that there is a proper-path connecting
$g$ and $g'$ in $G$, say $P=gg_1g_2,\cdots g_{\ell-1}g'$. Then the
path induced by the edges in
$\{(gh, g_1h), (g_1h, g_2h), \cdots (g_{\ell-1}h, g'h')\}$ is a
proper-path connecting $(g,h)$ and $(g',h')$. Therefore, the
above coloring is a proper-path coloring of $G\circ H$,
and hence $pc(G\circ H)=pc(G)=pc(H)$.

$(ii)$ If $pc(G)=1, pc(H)\geq 2$, then $pc(G\circ H)\leq 2$ by Lemma \ref{lem1} and Theorem \ref{th2}. Since $diam(G\circ H)\geq 2$, $pc(G\circ H)\geq 2$. So $pc(G\circ H)=2$

$(iii)$ The same as $(ii)$.

$(iv)$ If $pc(G)=1, pc(H)=1$, then both $G$ and $H$ are complete.
So $pc(G\circ H)=1$.
\qed
\end{pf}\\

To show the sharpness of the upper bound in Theorem \ref{th2}, we consider the following
example.\\

\noindent \textbf{Example 3:} Let $G=P_n$ be a path of order $n \ (n\geq 2)$ and $H=P_m$ be a path of order $m \ (m\geq 2)$.
If $m,n\geq 3$, then $pc(G)=pc(H)=2$, so
$pc(G\circ H)\leq pc(G)=pc(H)=2$ by Theorem \ref{th2}. Since $diam(G\circ H)\geq 2$, $pc(G\circ H)\geq 2$. So $pc(G\circ H)=2$.
If $m=2$, $n\geq 3$, then $pc(H)=1, pc(G)= 2$, so $pc(G\circ H)\leq 2$ by Theorem \ref{th2}. Since $diam(G\circ H)\geq 2$, it follows that $pc(G\circ H)\geq 2$. So $pc(G\circ H)=2$;
If $n=2$, $m\geq 3$, then $pc(G)=1, pc(H)= 2$, then $pc(G\circ H)=2 \leq 2$ by Theorem \ref{th2}. Since $diam(G\circ H)\geq 2$, we have $pc(G\circ H)\geq 2$. So $pc(G\circ H)=2$;
If $m=n=2$, then $pc(G)=1$, $pc(H)=1$ and $pc(G\circ H)=1$ by Theorem \ref{th2}.

\begin{cor}\label
{cor1}
Let $G$ and $H$ be connected graphs, then $pc(G\circ H)\leq \max\{pc(G), pc(H)\}$.
\end{cor}

\section{The direct product}

The direct product $G\times H$ of graphs $G$ and $H$ has the vertex
set $V(G)\times V(H)$. Two vertices $(g,h)$ and $(g',h')$ are
adjacent if the projections on both coordinates are adjacent, i.e.,
$gg'\in E(G)$ and $hh'\in E(H)$. It is clearly commutative and
associativity also follows quickly. For more general properties we
recommend \cite{Hammack}. The direct product is the most natural
graph product in the sense of categories. But this also seems to be
the reason that it is, in general, also the most elusive product of
all standard products. For example, $G\times H$ needs not to be
connected even when both factors are. To gain connectedness of
$G\times H$ at least one factor must additionally be nonbipartite as
shown by Weichsel \cite{Weichsel}. Also, the distance formula
$$
d_{G\times H}((g,h),(g',h')) = \min\{\max\{d^{e}_{G}(g,g'),d^{e}_{H}
(h,h')\},\max\{d^{o}_{G}(g,g'), d^{o}_{H}(h,h')\}\}
$$ for the direct
product is far more complicated as it is for other standard
products. Here $d^{e} _{G}(g,g')$ represents the length of a
shortest even walk between $g$ and $g'$ in $G$, and $d^{o}_{G}
(g,g')$ the length of a shortest odd walk between $g$ and $g'$ in
$G$. The formula was first shown in \cite{Kim} and later in
\cite{GhidewonH} in an equivalent version. There is no final
solution for the connectivity of the direct product, only some
partial results are known
 (see \cite{BresarS,GujiV}).

In this section we construct different upper bounds for the proper
connection number of the direct product with respect to some
invariants of the factors that are related to the rainbow
vertex-connection number of the factors. A similar concept as for
the distance formula is used and is due to the rainbow odd and even
walks between vertices (and not only rainbow paths) and is thus, in
a way, related with the formula.
We say that $G$ is \emph{odd-even proper connected} if there
exists a proper colored odd path and a proper
colored even path between every pair of (not necessarily different)
vertices of $G$. The \emph{odd-even proper connection
number} of a graph $G$, $oepv(G)$, is the smallest number of colors
needed for $G$ to be odd-even proper connected and it equals
infinity if no such a coloring exists. A bipartite graph has either
only even or only odd paths between two fixed vertices, thus there
is no odd-even proper coloring of such a graph. On the other
hand, let $G$ be a graph in which every vertex lies on some odd
cycle. Then $oepc(G)$ is finite since coloring every vertex with
its own color produces an odd-even proper coloring of $G$.

One can see that a odd cycle is an example where this coloring is
optimal, and $oervc(G)\leq |V(G)|$ for a connected graph $G$.

It is also easy to see that $oepc(K_{3})=3$. For $n\geq 3$, and $n$ is odd,
$oepc(C_{n})=3$. For $n\geq 3$, and $n$ is even,
$oepc(C_{n})=2$.

Let $G$ be a graph. We split $G$ into two spanning subgraphs $O^{G}$
and $B^{G}$, where the set $E(O^{G})$ consists of all edges of
$G$ that lie on some odd cycle of $G$, and the set $E(B^{G}) =
E(G)\setminus E(O^{G})$. Clearly, $O^{G}$ and $B^{G}$ are not always
connected. Let $O^{G}_{ 1} ,O^{G}_{ 2 },\cdots,O^{G}_{ k}$ and
$B^{G}_{1},B^{G}_{2},\cdots,B^{G}_{\ell}$ be components of $O^{G}$
and $B^{G}$, respectively, each one containing more than one vertex.
Let
$$
o(G)=oepc(O^{G}_{1})+oepc(O^{G}_{2})+\cdots+
oepc(O^{G}_{k}),
$$
and
$$
b(G)=pc(B^{G}_{1})+pc(B^{G}_{2})+\cdots+
pc(B^{G}_{\ell})
$$

Note that $o(G)$ is finite since it is defined on nontrivial
components $O^{G}_{i}$, $i\in \{1,2,\cdots,k\}$.

\begin{thm}\label{th3}
Let $G$ and $H$ be a nonbipartite connected graph. Then
$$
pc(G\times H)\leq \min\{pc(H)((b(G)+o(G)), pc(G)(b(H)+o(H))\}.
$$
\end{thm}
\begin{pf}
Without loss of generality, $pc(H)((b(G)+o(G))\leq pc(G)(b(H)+o(H))$.
Denote by $c^{B}_{G}$ an optimal proper-coloring of
components of $B^{G}$. Let $c^{O}_{G}$ be an optimal odd-even
proper-coloring of components of $O^{G}$.

We give a proper-coloring of $G\times H$ as follows.
If $e\in E(G\times H)$ projects on
$G$ to $e'\in B_G$, we set $c(e)=(c^{B}
_{G}(e'),c_H(e''))$,
and if $e$ projects on $G$ to $e'\in O_G$,
we set $c(e)=(c^{O}_{G}(e'),c_H(e''))$. where $e''\in E(H)$ is the projection of $e$ on $H$.
By this way, we get a coloring of $V(G\times H)$ with
$pc(H)(o(G)+b(G))$ colors and it remains to show that this is a
rainbow coloring of $G\times H$.

Let $(g, h)$ and $(g', h')$ be arbitrary vertices from $G\times H$.
Clearly, there is a proper path connecting $g$ and $g'$, say
$P=gg_{1},\ldots g_{\ell-1}g'$. By the same reason, there is a
proper path connecting $h$ and $h'$, say $Q=hh_{1},\ldots
h_{k-1}h'$. Observe that $P$ is a shortest proper $g,
g'$-path in $G$ induced by ${B}_{G}$ and ${O}_{G}$, and $Q$ is a
shortest proper $h, h'$-path in $H$. If $g=g'$ or $h=h'$,
then $P$ or $Q$, respectively, is a trivial one vertex path.

We distinguish the following two cases to prove this theorem.

\emph{Case $1$.} $\ell$ and $k$ have the same parity.

If $h=h'$, then we let $h_{k-1}$ be an arbitrary neighbor of $h$.
Then the path induced by the edges in
$$
\{(gh, g_{1}h_{k-1}),(g_{1}h_{k-1}, g_{2}h),(g_{2}h,g_{3}
h_{k-1}),\ldots,(g_{\ell-1}h_{k-1}, g'h')\}
$$
is a proper$(g, h), (g', h')$-path in $G\times H$.

If $g=g'$, then we let $g_{\ell-1}$ be an arbitrary neighbor of $g$.
Then the path induced by the edges in
$$
\{(gh, g_{\ell-1}h_{1}),(g_{\ell-1}h_{1}, gh_{2}),(gh_{2}, g_{\ell-1},
h_{3}),\ldots,(g_{\ell-1}h_{k-1}, g'h')\}
$$
is a vertex-rainbow $(g, h), (g', h')$-path in $G\times H$.

If $g\neq g'$, and $h\neq h'$, then the path induced by the edges in
$$
\{(gh, g_{1}h_{1}), (g_{1}h_{1}, g_{2}h_{2})\ldots, (g_{k}h', g_{k+1}h_{k-1}),(g_{k+1}h_{k-1},g_{k+2}h')\ldots
 (g_{\ell-1}h_{k-1}, g'h')\}
$$
is a proper $(g, h), (g', h')$-path in $G\times H$ whenever
$\ell\geq k$, and the path induced by the edges in
$$
\{(gh, g_{1}h_{1}), (g_{1}h_{1}, g_{2}h_{2})\ldots, (g_{\ell-1}h_{\ell-1}, g'h_{\ell}),
(g'h_{\ell}, g_{\ell-1}h_{\ell+1})\ldots, (g_{\ell-1}h_{k-1}, g'h')\}
$$
is a proper $(g, h), (g', h')$-path in $G\times H$ whenever
$\ell< k$.

\emph{Case 2.} $\ell$ and $k$ have different parity.

If there exists a $g_{i}, g_{j}$-subpath of $P$ in $O^{G}_{p}$, we
replace this subpath by a rainbow $g_{i}, g_{j}$-path of different
parity in $O^{G}_{p}$ to obtain a proper path $P'$ between $g$ and
$g'$. If this is the case, then $|E(P')|$ and $k$ have the same
parity and we can use Case $1$.
We now assume that all the $g_{i}, g_{j}$-subpaths of $P$ in
$B^{G}_{p}$, that is, all vertices of $P$ are in $B^{G}_{p}$. To
find a proper $(g, h), (g', h')$-path in $G\times H$, we
find out a $g,g'$-walk in $G$. Note that $P$ is contained in one
component $B^{G}_{q}$.
Let $g_{i}\in V(P)$ be a vertex that is closest to any component
$O^{G}_{p}$ of $G$ and let $v_{1}\in O^{G}_{p}$ be closest to
$g_{i}$. Let $R = g_{i}g'_{i+1},\ldots, g'_{i+r} \ (g'_{i+r}=v_{1})$
be a shortest $g_{i}, v_{1}$-path. From the definition of odd-even
rainbow vertex-coloring, we know that there exists an odd
vertex-rainbow $v_{1}, v_{1}$-cycle $C = v_{1}v_{2},\ldots
v_{p}v_{1}$ in $O^{G}_{p}$. Now we insert a closed walk that follows
$RCR$ from $g_{i}$ into a path $P$ to obtain a $g, g'$-walk
\begin{eqnarray*}
W&=&gg_{1}\ldots g_{i}g'_{i+1}\ldots g'_{i+r}v_{2}v_{3},\ldots
v_{p}v_{1}g'_{i+r-1}g'_{i+r-2}\ldots g'_{i+1}g_{i}g_{i+1}\ldots
g'\\
&=&u_{0}u_{1},\ldots u_{\ell+p+2r}.
\end{eqnarray*}
of length $t=\ell+2r+p$. Note that $t$ and $\ell$ have different
parity since $p$ is an odd number, and thus $t$ and $k$ have the
same parity. If $k\geq t$, then the path induced by the edges in
$$
\{(u_0h, u_1h_1), (u_1h_1, u_2h_2),\cdots (u_th_t, u_{t-1}h_{t+1}), (u_{t-1}h_{t+1}, u_th_{t+2}),\cdots (u_{t-1}h_{k-1}, u_th')\}
$$
is a proper-coloring connected $gh$ and $g'h'$.

If $k<t$, then the path induced by the edges in
$$
\{(u_0h, u_1h_1), (u_1h_1, u_2h_2),\cdots (u_{k-1}h_{k-1}, u_kh'), (u_kh', u_{k+1}h_{k-1}),\cdots (u_{t-1}h_{k-1}, u_th')\}
$$
is a proper-coloring connected $gh$ and $g'h'$.
\qed
\end{pf}

\begin{cor}\label{cor2}
Let $G$ and $H$ be connected graphs, where $G$ is nonbipartite
and $H$ is bipartite. Then
$$
pc(G\times H)\leq pc(H)(b(G)+o(G)).
$$
\end{cor}

A bipartite graph $G=(V_{0}\cup V_{1},E)$ is said to \emph{have a
property $\pi$} if $G$ admits of an automorphism $\psi$ such that
$x\in V_{0}$ if and only if $\psi(x)\in V_{1}$. For more details, we
refer to \cite{JhaKZ}.

\begin{lem}\label{lem4}{\upshape \cite{JhaKZ}}
If $G$ and $H$ are bipartite graphs one of which has property $\pi$,
then the two components of $G\times H$ are isomorphic.
\end{lem}

\begin{pro}\label{pro3}
Let $G$ be a nonbipartite connected graph. Then 
$$
pc(G\times
K_{2})\leq o(G)+b(G).
$$
\end{pro}
\begin{pf}
 Let $c^{O}_{G}$ be an optimal odd-even proper-coloring of $O^{_{G}}$ and let $c^{B}_{G}$ be an optimal proper-coloring of $B^{_{G}}$ (for both cases it holds that no color
appears in two different components). Observe that $c^{O}_{G}=o(G)$
and $c^{B}_{G}=b(G)$. We provide a coloring $c$ of $G\times
K_{2}$ with $o(G)+b(G)$ colors as follows.

Recall that $O^{G}_{1} ,O^{G}_{2},\cdots,O^{G}_{k}$ and
$B^{G}_{1},B^{G}_{2},\cdots,B^{G}_{\ell}$ are all the components of
$O^{G}$ and $B^{G}$, respectively. By the definition, $B^{G}_{i}$ is
bipartite graph. From Lemma \ref{lem3}, $B^{G}_{i}\times K_{2}$ can
be decomposed into two subgraphs isomorphic to $B^{G}_{i}$. Color
both components of $B^{G}_{i}\times K_{2}$ (which are isomorphic to
$B^{G}_{i}$) optimally with $pc(B^{G}_{i})$ colors for every $i\in
\{1,2,\ldots, \ell\}$. For this we use $b(G)$ colors. Now, we
assign $o(G)$ new colors to the remaining vertices. For an edge
$(gh, g'h')$ of $G\times K_{2}$, it project on $G$ to an edge $gg'$ of
$O^{G}$ receive color $c(gh, g'h')=c^{O}_{G}(gg')$. For an edge
$(gh, g'h')$ of $G\times K_{2}$, it project on $G$ to an edge $gg'$ of
$B^{G}$ receive color $c(gh, g'h')=c^{B}_{G}(gg')$.
For the introduced
coloring $o(G)+b(G)$ colors are used and we need to show that $c$
is a proper-coloring of $G\times K_{2}$.

Set $V(K_{2})= \{k_1,k_2\}$. Let $(g,h)$ and $(g',h')$ be arbitrary
vertices in $G\times K_{2}$. Let $P=gg_{1},\ldots g_{\ell-1}g'$ be a
proper $g,g'$-path under the proper-coloring of $G$
induced by $c^{O}_{G}$ and $c^{B}_{G}$. We distinguish two cases to
show this proposition.

\emph{Case $1$.} Let $\ell$ and $d_{K_{2}}(h,h')$ have the same
parity.

Without loss of generality we may assume that $h=k_{1}$.
Consequently $h'= k_{1}$ if $\ell$ is an even number and $h'=k_{2}$
otherwise. Thus
$$
(gk_{1})(g_{1}k_{2})(g_{2}k_{1})\cdots(g'h')
$$
is a proper $(g,h),(g',h')$-path in $G\times K_{2}$.

\emph{Case $2$.} Let $\ell$ and $d_{K_{2}}(h,h')$ have different
parity.

Suppose first that $P$ has a nonempty intersection with some $O^{G}
_{p}$ and let $g_{i}$ be the first and $g_{j}$ the last vertex of
$P$ in $O^{G}_{p}$. Then we can find a proper
$g_{i},g_{j}$-path in $O^{G}_{p}$ with length of different parity as
is the length of the $g_{i},g_{j}$-subpath of $P$ in $O^{G}_{p}$.
Replacing the $g_{i},g_{j}$-subpath of $P$ by this proper
$g_{i},g_{j}$-path in $O^{G}_{p}$ we obtain a proper
$g,g'$-path of the same parity as $d_{K_{2}}(h,h')$ and we continue
as in Case $1$.

Suppose now that $P$ has an empty intersection with every
$O^{G}_{p}$, $p\in\{1,2,\ldots,k\}$. Then $P$ is contained in
$B^{G}_{q}$ for some $q$, and $(g,h)$ and $(g',h')$ are in different
components $(B^{G}_{ q})_{1}$ and $(B^{G}_{q})_{2}$ of
$B^{G}_{q}\times K_{2}$, respectively. Since $G$ is nonbipartite,
there exists a vertex $g''$ in some component of $O^{G}_{p}$. Set
$\{h_{r},h_{s}\}=\{k_{1},k_{2}\}$. Take a proper path from $(g,h)$
to $(g'',h_{r})$ in $(B^{G}_{q} )_{1}$, a proper odd path from
$(g'',h_{r})$ to $(g'',h_{s})$ in $O^{G}_{p}$ , and a rainbow path
from $(g'',h_{s})$ to $(g',h')$ in $(B^{G}_{q})_{2}$. This is a
proper $(g,h),(g',h')$-path in $G\times K_{2}$ since we have
used different colors for $(B^{G}_{q} )_{1}, (B^{G}_{q })_{2}$, and
$O^{G}_{ p}$.\qed
\end{pf}

\section{Applications}

In this section, we demonstrate the usefulness of the proposed
constructions by applying them to some instances of Cartesian and
lexicographical product networks.

The following results will be used later.

\begin{lem}{\upshape \cite{Hammack}}\label{lem4}
Let $(gh)$ and $(g'h')$ be two vertices of $G\circ H$.
Let $d_{G}(g)$ denote the degree of vertex $g$ in
$G$. Then
$$
d_{G\circ H}(gh,g'h')=\left\{
\begin{array}{ll}
d_{G}(gg'), &if~g\neq g';\\[0.2cm]
d_{H}(hh'), &if~g=g'~and~d_{G}(g)=0;\\[0.2cm]
\min \{d_{H}(hh'),2\}, &if~g=g'~and~d_{G}(g)\neq 0.
\end{array}
\right.
$$
\end{lem}

\subsection{Two-dimensional grid graph}

A \emph{two-dimensional grid graph} is an $m\times n$ graph
$G_{n,m}$ that is the graph Cartesian product $P_n\Box P_m$ of path
graphs on $m$ and $n$ vertices. See Figure 1 $(a)$ for the case
$m=3$. For more details on grid graph, we refer to \cite{Calkin,
Itai}. The network $P_n\circ P_m$ is the graph lexicographical
product $P_n\circ P_m$ of path graphs on $m$ and $n$ vertices. For
more details on $P_n\circ P_m$, we refer to \cite{Mao2}. See Figure
1 $(b)$ for the case $m=3$.

\begin{figure}[!hbpt]
\begin{center}
\includegraphics[scale=0.9]{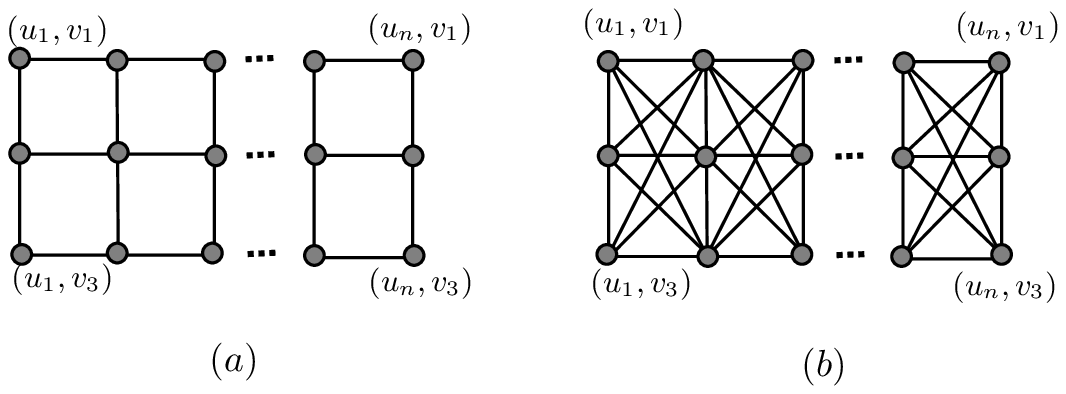}\\
Figure 1: $(a)$ Two-dimensional grid graph $G_{n,3}$; $(b)$ The
network $P_n\circ P_3$.
\end{center}\label{fig7}
\end{figure}

\begin{pro}\label{pro3}
$(i)$ For network $P_n\Box P_m \ (n\geq 2, m\geq 2)$,
$2\leq pc(P_n\Box P_m)\leq 3$.

$(ii)$ For network $P_n\circ P_m $,
$pc(P_n\circ P_m)=1$ when $m=n=2$, $pc(P_n\circ P_m)=2$ when $m=2, n>2$ or $n=2, m>2$ or $m,n>2$.
\end{pro}
\begin{pf}
$(i)$ By Theorem \ref{th1}, we have $pc(P_n\Box P_m)\leq \min\{pc(P_n), pc(P_m)\}+1=2+1=3$. Observe that $diam(P_n\Box P_m)\geq 2$. So $2\leq pc(P_n\Box P_m)\leq 3$.

$(ii)$ The same as Example $3$. \qed
\end{pf}

\subsection{$n$-dimensional mesh}

An \emph{$n$-dimensional mesh} is the Cartesian product of $n$
linear arrays. By this definition, two-dimensional grid graph is a
$2$-dimensional mesh. An $n$-dimensional hypercube is a special case
of an $n$-dimensional mesh, in which the $n$ linear arrays are all
of size $2$; see \cite{Johnsson}.

\begin{pro}\label{pro4}
$(i)$ For $n$-dimensional mesh $P_{L_1}\Box P_{L_2}\Box \cdots \Box
P_{L_n}$,
$$
pc(P_{L_1}\Box P_{L_2}\Box \cdots \Box
P_{L_n})=2.
$$

$(ii)$ For network $P_{L_1}\circ P_{L_2}\circ
\cdots \circ P_{L_n}$,
if there exists some $L_j$ such that $L_j\neq 2 \ (1\leq j\leq n)$, then $pc(P_{L_1}\circ P_{L_2}\circ
\cdots \circ P_{L_n})=2$;
If $L_1=L_2=\cdots =L_n=2$, then $pc(P_{L_1}\circ P_{L_2}\circ
\cdots \circ P_{L_n})=1$.
\end{pro}
\begin{pf}
$(i)$ By Lemma \ref{lem1}, we have $diam((P_{L_1}\Box P_{L_2}\Box
\cdots \Box P_{L_n})=\sum_{i=1}^ndiam(P_{L_i})=\sum_{i=1}^n(L_i-1)$
$=\sum_{i=1}^nL_i-n\geq 2$. By Theorem \ref{th2},
$pc(P_{L_1}\Box P_{L_2}\Box \cdots \Box
P_{L_n})\leq \min\{pc(P_{L_1})$, $pc(P_{L_2}), \cdots, pc(P_{L_n})\}+1=2$.
So $pc(P_{L_1}\Box P_{L_2}\Box \cdots \Box
P_{L_n})=2.$

$(ii)$ If there exists some $L_j$ such that $L_j\neq 2 \ (1\leq j\leq n)$, then $pc(P_{L_1}\circ P_{L_2}\circ
\cdots \circ P_{L_n})\leq \max\{P_{L_1}, P_{L_2},\cdots P_{L_n}\}=2$ by Corollary\ref{cor1}.
Since $diam(P_{L_1}\circ P_{L_2}\circ
\cdots \circ P_{L_n})\geq 2$, $pc(P_{L_1}\circ P_{L_2}\circ
\cdots \circ P_{L_n})\leq 2$. So $pc(P_{L_1}\circ P_{L_2}\circ
\cdots \circ P_{L_n})=2$.

If $L_1=L_2=\cdots =L_n=2$, then $P_{L_1}\circ P_{L_2}\circ
\cdots \circ P_{L_n}$ is a complete graph. So $pc(P_{L_1}\circ P_{L_2}\circ
\cdots \circ P_{L_n})=1$.
 \qed
\end{pf}

\subsection{$n$-dimensional torus}

An \emph{$n$-dimensional torus} is the Cartesian product of $n$
rings $R_1,R_2,\cdots,R_n$ of size at least three.(A
ring is a cycle in Graph Theory.) The rings $R_i$ are not necessary
to have the same size. Ku et al.
\cite{Ku} showed that there are $n$ edge-disjoint spanning trees in
an $n$-dimensional torus. The network $R_1\circ R_2\circ \cdots
\circ R_n$ is investigated in \cite{Mao2}. Here, we consider the
networks constructed by $R_1\Box R_2\Box \cdots \Box R_n$ and
$R_1\circ R_2\circ \cdots \circ R_n$.

\begin{pro}\label{pro7}
$(i)$ For network $R_1\Box R_2\Box \cdots \Box R_n$,
$$
2\leq pc(R_1\Box R_2\Box \cdots
\Box R_n)\leq \min\{pc(R_1), pc(R_2), \cdots pc(R_n)\}+1=3
$$
where $r_i$ is the order of $R_i$ and $3\leq i\leq
n$.

$(ii)$ For network $R_1\circ R_2\circ \cdots \circ R_n$,
$$
pc(R_1\circ R_2\circ \cdots
\circ R_n)=2.
$$
\end{pro}
\begin{pf}
$(i)$ By Lemma \ref{lem1}, we have $diam(R_1\Box R_2\Box \cdots \Box
R_n)=\sum_{i=1}^ndiam(R_i)=\sum_{i=1}^n\lfloor r_i/2\rfloor\geq 2$ and hence $pc(R_1\Box R_2\Box \cdots \Box R_n)\geq 2$. By Theorem \ref{th1}, we have
$$
pc(R_1\Box R_2\Box \cdots
\Box R_n)\leq \min\{pc(R_1), pc(R_2), \cdots pc(R_n)\}+1=3.
$$
Therefore, $2\leq pc(R_1\Box R_2\Box \cdots
\Box R_n)\leq 3$.

$(ii)$ From Corollary\ref{cor1}, we have
$pc(R_1\circ R_2\circ \cdots \circ R_n)\leq \max \{pc(R_1), pc(R_2),\cdots pc(R_n)\}=2$. Since $diam(R_1\circ R_2\circ \cdots \circ R_n)\geq 2$, $pc(R_1\circ R_2\circ \cdots \circ R_n)\geq 2$. So $pc(R_1\circ R_2\circ \cdots \circ R_n)=2$.\qed
\end{pf}

\subsection{$n$-dimensional generalized hypercube}

Let $K_m$ be a clique of $m$ vertices, $m\geq 2$. An
\emph{$n$-dimensional generalized hypercube} \cite{DayA,
Fragopoulou} is the Cartesian product of $m$ cliques. We have the following:

\begin{pro}\label{pro6}
$(i)$ For network $K_{m_1}\Box K_{m_2}\Box \cdots \Box K_{m_n} \
(m_i\geq 2, \ n\geq 2, \ 1\leq i\leq n)$
$$
pc(K_{m_1}\Box K_{m_2}\Box \cdots \Box K_{m_n})=2
$$

$(ii)$ For network $K_{m_1}\circ K_{m_2}\circ \cdots \circ K_{m_n}$,
$$
pc(K_{m_1}\circ K_{m_2}\circ \cdots \circ K_{m_n})=1.
$$
\end{pro}
\begin{pf}
$(1)$ Observe that $diam(K_{m_1}\Box K_{m_2}\Box \cdots \Box
K_{m_n})=\sum_{i=1}^ndiam(K_{m_i})=n\geq 2$. So $pc(K_{m_1}\Box K_{m_2}\Box \cdots \Box
K_{m_n})\geq 2$. By Theorem \ref{th1}, we have $pc(K_{m_1}\Box K_{m_2}\Box \cdots \Box K_{m_n})\leq \min\{pc(K_{m_1}), pc(K_{m_2})\cdots pc(K_{m_n})\}+1=2$. So $pc(K_{m_1}\Box K_{m_2}\Box \cdots \Box K_{m_n})=2$.

$(2)$ Observe that $K_{m_1}\circ K_{m_2}\circ \cdots
\circ K_{m_n}$ is a complete graph. So $pc(K_{m_1}\circ K_{m_2}\circ \cdots \circ
K_{m_n})=1$.\qed
\end{pf}

\subsection{$n$-dimensional hyper Petersen network}

An \emph{$n$-dimensional hyper Petersen network} $HP_n$ is the
Cartesian product of $Q_{n-3}$ and the well-known Petersen graph
\cite{Das}, where $n\geq 3$ and $Q_{n-3}$ denotes an
$(n-3)$-dimensional hypercube. The cases $n=3$ and $4$ of hyper
Petersen networks are depicted in Figure 5. Note that $HP_3$ is just
the Petersen graph (see Figure 5 $(a)$).

The network $HL_n$ is the lexicographical product of $Q_{n-3}$ and
the Petersen graph, where $n\geq 3$ and $Q_{n-3}$ denotes an
$(n-3)$-dimensional hypercube; see \cite{Mao2}. Note that $HL_3$ is
just the Petersen graph, and $HL_4$ is a graph obtained from two
copies of the Petersen graph by add one edge between one vertex in a
copy of the Petersen graph and one vertex in another copy. See
Figure 5 $(c)$ for an example (We only show the edges $v_1u_i \
(1\leq i\leq 10)$).

\begin{figure}[!hbpt]
\begin{center}
\includegraphics[scale=0.8]{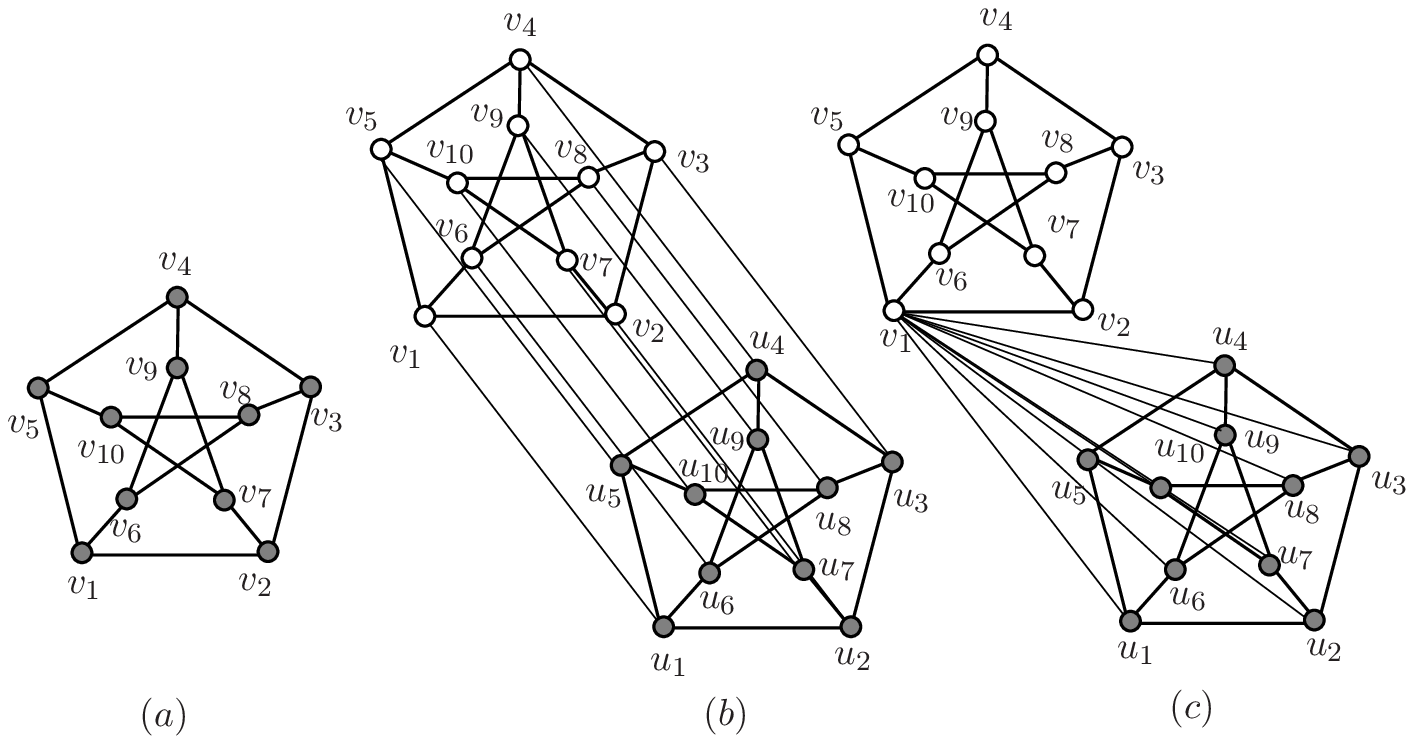}\\
Figure 2: $(a)$ Petersen graph; $(b)$ The network $HP_4$; $(c)$ The
structure of $HL_4$.
\end{center}\label{fig7}
\end{figure}

\begin{pro}\label{pro7}
$(1)$ For network $HP_3$ and $HL_3$, $pc(HP_3)=pc(HL_3)=2$;

$(2)$ For network $HL_4$ and $HP_4$, $2\leq pc(HP_4)\leq 3$ and $pc(HL_4)= 2$.
\end{pro}
\begin{pf}
$(1)$ Since $diam(HP_3)=diam(HL_3)=2$, it follows that $pc(HP_3)=pc(HP_3)\geq 2$.
One can check that there is a proper-coloring with two colors. So $pc(HP_3)=pc(HL_3)=2$.

$(2)$ From Theorem \ref{th1}, $pc(HP_4)\leq 3$. Since $diam(HP_4)=2$, it follows that $pc(HP_4)\geq 2$. So $pc(HP_4)= 2$.
From Corollary \ref{cor1}, we have $pc(HL_4)\leq 2$. Since $diam(HL_4)=2$, $pc(HL_4)\geq 2$. So $pc(HL_4)= 2$.\qed
\end{pf}


\begin{thebibliography}{1}

\bibitem{AndrewsLLZ}
E. Andrews, E. Laforge, C. Lumduanhom, P. Zhang, \emph{On proper-path colorings in
graphs}, J. Combin. Math. Combin. Comput, to appear.

\bibitem{AnandCKP}
B.S. Anand, M. Changat, S. Klav\u{z}ar, I. Peterin, \emph{Convex
sets in lexicographic products of graphs}, Graphs Combin. 28(2012),
77--84.

\bibitem{Bondy} J.A. Bondy, U.S.R. Murty,
{\it Graph Theory}, GTM 244, Springer, 2008.

\bibitem{BresarS}
B. Bre\v{s}ar, S. \v{S}pacapan, \emph{On the connectivity of the
direct product of graphs}, Australas. J. Combin. 41(2008), 45--56.


\bibitem{Calkin}
N.J. Calkin, H.S. Wilf, \emph{The number of independent sets in a
grid graph }, SIAM J. Discrete Math. 11(1)(1998), 54--60.

\bibitem{CFMY}
S. Chakraborty, E. Fischer, A. Matsliah, R. Yuster, \emph{Hardness
and algorithms for rainbow connectivity}, 26th International
Symposium on Theoretical Aspects of Computer Science STACS (2009),
243--254. Also, see {\it J. Combin. Optim.} 21(2011), 330--347.

\bibitem{Chartrand}
G. Chartrand, G.L. Johns, K.A. McKeon, P. Zhang, \emph{Rainbow
connection in graphs}, Math. Bohem. 133(2008), 85--98.

\bibitem{ChartrandJMZ}
G. Chartrand, G. L. Johns, K. A. McKeon, P. Zhang, \emph{The rainbow
connectivity of a graph}, Networks 54(2009), 75--81.

\bibitem{ChartrandOZ}
G. Chartrand, F. Okamoto, P. Zhang, \emph{Rainbow trees in graphs and
generalized connectivity}, Networks 55(2010), 360--367.

\bibitem{Das}
S.K. Das, S.R. \"{O}hring, A.K. Banerjee, \emph{Embeddings into
hyper Petersen network: Yet another hypercube-like interconnection
topology}, VLSI Design, 2(4)(1995), 335--351.

\bibitem{DayA}
K. Day, A.-E. Al-Ayyoub, \emph{The cross product of interconnection
networks}, IEEE Trans. Parallel and Distributed Systems 8(2)(1997),
109--118.

\bibitem{Fragopoulou}
P. Fragopoulou, S.G. Akl, H. Meijer, \emph{Optimal communication
primitives on the generalized hypercube network}, IEEE Trans.
Parallel and Distributed Computing 32(2)(1996), 173--187.

\bibitem{Hammack}
R. Hammack, W. Imrich, Sandi Klav\u{z}r, \emph{Handbook of product
graphs}, Secend edition, CRC Press, 2011.


\bibitem{M.Krivelevich}
M. Krivelevich, R. Yuster,
\emph{The rainbow connection of a graph is (at~most) reciprocal to
its minimum degree three}, IWOCA 2009, LNCS 5874(2009), 432--437.



\bibitem{LiSS}
X. Li, Y. Shi, Y. Sun, \emph{Rainbow connections of graphs--A survey},
Graphs Combin. 29(1)(2013), 1--38.


\bibitem{LiS}
X. Li, Y. Sun, {\it Rainbow Connections of Graphs}, SpringerBriefs
in Math., Springer, New York, 2012.

\bibitem{LiWY}
X. Li, M. Wei, J. Yue, \emph{Proper connection number and connected dominating sets},
arXiv 1501. 05717 v1 [math. CO] 23 Jan 2015.


\bibitem{Gologranc}
T. Gologranc, Ga\v{s}per Meki\v{s}, I. Peterin, \emph{Rainbow
connection and graph products}, 30(3)(2014), 591--607.


\bibitem{GhidewonH}
A.A. Ghidewon, R. Hammack, \emph{Centers of tensor product of
graphs}, Ars Combin. 74(2005), 201--211.

\bibitem{GujiV}
R. Guji, E. Vumar, \emph{A note on the connectivity of Kronecker
products of graphs}, Appl. Math. Lett. 22(2009), 1360--1363.

\bibitem{HuangLW}
F. Huang, X. Li, S. Wang, \emph{Proper connection numbers of
complementary graphs}, arXiv 1504. 02414 v2 [math. CO] 29 Apr 2015.

\bibitem{Itai}
A. Itai, M. Rodeh, \emph{The multi-tree approach to reliability in
distributed networks}, Information and Computation 79(1988), 43--59.


\bibitem{JhaKZ}
P.K. Jha, S. Klav\v{z}ar, B. Zmazek, \emph{Isomorphic components of
Kronecker product of bipartite graphs}, Discuss. Math. Graph Theory
17(1997), 301--309.

\bibitem{Johnsson}
S.L. Johnsson, C.T. Ho, \emph{Optimum broadcasting and personaized
communication in hypercubes}, IEEE Trans. Computers 38(9)(1989),
1249--1268.

\bibitem{Kim}
S.R. Kim, \emph{Centers of a tensor composite graph}, Congr. Numer.
81 (1991) 193--203.

\bibitem{KlavzarM}
S. Klav\v{z}ar, G. Meki\v{s}, \emph{On the rainbow connection of
Cartesian products and their subgraphs}, Discuss. Math. Graph Theory
32 (2012), 783--793.

\bibitem{Klavzar}
S. Klav\v{z}ar, S. \v{S}pacapan, \emph{On the edge-connectivity of
Cartesian product graphs}, Asian-Eur. J. Math. 1 (2008), 93--98.

\bibitem{Krivelevich}
 M. Krivelevich, R. Yuster, \emph{The rainbow connection of a graph is
(at most) reciprocal to its minimum degree}, J. Graph Theory 63
(2009), 185--191.

\bibitem{Ku} S. Ku, B. Wang, T. Hung, \emph{Constructing edge-disjoint
spanning trees in product networks}, Parallel and Distributed
Systems, IEEE Transactions on parallel and disjoited systems 14(3)
(2003), 213-221.

\bibitem{Mao2}
Y. Mao, \emph{Path-connectivity of lexicographical product  graphs},
Int. J. Comput. Math., in press.

\bibitem{Nowakowski}
R.J. Nowakowski, K. Seyffarth, \emph{Small cycle double covers of
products. I. Lexicographic product with paths and cycles}, J. Graph
Theory 57 (2008), 99--123.


\bibitem{Spacapan}
S. \v{S}pacapan, \emph{Connectivity of strong products of graphs},
Graphs Combin. 26 (2010), 457--467.

\bibitem{Weichsel}
P.M. Weichsel, \emph{The Kronecker product of graphs}, Proc. Amer.
Math. Soc. 13(1962), 47--52.

\bibitem{Zhu}
X. Zhu, \emph{Game coloring the Cartesian product of graphs}, J.
Graph Theory 59(2008), 261--278.













\end{thebibliography}
\end{document}